\documentclass[a4paper,10pt]{amsart}
\usepackage{standalone}
\usepackage{amsfonts}
\usepackage{amsthm}
\usepackage{amssymb}
\usepackage{amsmath}
\usepackage{theoremref}
\usepackage{enumerate}
\usepackage{bbm}
\usepackage{bm}
\usepackage{pgfplots}
\pgfplotsset{compat=1.18} 
\usepgfplotslibrary{fillbetween}
\usetikzlibrary{intersections}
\usepackage{mathtools}
\usepackage{tkz-euclide}
\usetikzlibrary{patterns}
\usepackage{hyperref}
\usepackage{tikz}
\usetikzlibrary{decorations.pathreplacing}
\usepackage[skins,breakable]{tcolorbox}
\usepackage[
    a4paper,
    left=2cm,
    right=2cm,
    top=1.7cm,
    bottom=1.7cm
]{geometry}

\newcommand{\T}{\mathbb{T}}

\newcommand{\D}{\mathbb{D}}

\newcommand{\dist}[2]{\text{dist}( #1, #2 ) }

\newcommand{\m}{\textit{m}}

\renewcommand{\k}{\mathcal{K}}

\newcommand{\Carl}{\operatorname{Carl}}

\newtheorem{thm}{Theorem}
\newtheorem*{thm*}{Theorem}

\newtheorem{cor}[thm]{Corollary}
\newtheorem*{cor*}{Corollary}
\newtheorem{prop}[thm]{Proposition}
\theoremstyle{definition}

\theoremstyle{definition}

\newtheorem*{defn*}{Definition}

\newtheorem{claim*}{Claim}

\title{Polar Carleson zero set preventing Dirichlet cyclicity}
\author{Bartosz Malman}

\address{Institutionen för ekonomi och matematik, 
        Mälardalens Universitet,
		721 23 Västerås, Sweden}
\email{bartosz.malman@mdu.se} 

\begin{document}

\begin{abstract}
    The work \cite{chen2026brownshields} of Chen, Fang and Guo shows, finally, that the Brown-Shields conjecture is false. In this note\footnotemark[1], we  construct a polar Carleson set $E$ on the unit circle $\T$ which has the property that any Hölder continuous function in the Dirichlet space vanishing on $E$ is non-cyclic.
\end{abstract}

\maketitle
\setcounter{footnote}{1}

\footnotetext[1]{All substantive results in this note have been given by the Oracle, in this case the ChatGPT 6.0 Astra model. But she did not write the note.}

\subsection{Reduction to an estimate of a point evaluation and a discussion}

Let $\mathcal{D}$ be the Dirichlet space on the unit disk $\D$, equipped with the norm $\| f \|_{\mathcal{D}}^2 = |f(0)|^2 + \int_\D |f'(z)|^2\,dA(z)$, where $dA$ is normalized area measure. On $\T$, let $m$ be normalized Lebesgue measure and write $|B|=m(B)$ for every measurable set $B\subseteq\T$. 

\setcounter{equation}{1}
\begin{prop}\thlabel{P:BScounterexample}
    There exist a compact Carleson set $E\subset\T$ of zero logarithmic capacity and a numerical constant $C>0$ with the following property. If $h\in C(\T) \cap\mathcal D$ vanishes on $E$ and satisfies the distance estimate \[|h(t)| \leq C_h \dist{t}{E}^{1/2}, \quad t \in \T,\] then
    \begin{equation}\label{E:SmallAverageFromVanishing}
        |h(0)|^2 \leq C\int_\D |h'(z)|^2\,dA(z).
    \end{equation}
    In particular, the estimate in \eqref{E:SmallAverageFromVanishing} does not involve the constant $C_h$.\footnote{We mention, at once, that the magical disappearance of dependence on $C_h$ occurs in the limit as $p \to \infty$ in \eqref{E:IterationEq} in which an iterative inequality associated to $E$ is employed.}
\end{prop}

The definition of Carleson sets is given below, right after \eqref{E:CarlesonEntropyDef}. We spend the note on proving \thref{P:BScounterexample}.

For $\alpha \in (0,1)$, let $C^{\alpha}(\T)$ be the space of Hölder regular functions on $\T$ with exponent $\alpha$, namely functions satisfying an estimate of the form $|f(t)-f(s)| \leq C_{f,\alpha} |t-s|^\alpha$ for some constant $C_{f,\alpha} > 0$. Then the distance estimate in \thref{P:BScounterexample} is satisfied for $h \in C^{1/2}(\T)$ that vanishes on $E$. If $h \in C^{\alpha}(\T)$ vanishes on $E$ and the integer $k$ satisfies $k\alpha \geq 1/2$, then the distance estimate is satisfied for $h^k$ and $p h^k$, where $p$ is any polynomial.

\begin{cor}\thlabel{C:Noncyclicity}
Every outer function $f\in C^{\alpha}(\T)\cap\mathcal D$ that vanishes on the set $E$ in \thref{P:BScounterexample} is non-cyclic in $\mathcal D$.
\end{cor}

\begin{proof}
Let $k$ be an integer satisfying $k\alpha \geq 1/2$, and set $g = f^k$. Since $f$ is bounded, we have $g \in \mathcal{D}$. By the paragraph preceding the corollary, we have
\[ |p(0)g(0)|^2 \leq C\int_\D |(p(z)g(z))'|^2 dA(z)\] for any polynomial $p$. Since
\[ \int_\D |(p(z)g(z))'|^2 dA(z) = \int_\D |(p(z)g(z)-1)'|^2 dA(z) \leq \|pg-1\|^2_{\mathcal{D}}\] we deduce that
\[ |p(0)g(0)|^2 \leq C \|pg - 1\|^2_{\mathcal{D}}\] and we see that we cannot have $pg \to 1$ in $\mathcal{D}$, for then the left-hand side above tends to $1$ while the right hand side tends to $0$.
Hence $g$ is non-cyclic. By \cite[Theorem 9.1.7]{dirichletspaceprimer}, $g = f^k$ is cyclic if and only if $f$ is cyclic.
\end{proof}

Given a Carleson set $E$ with $|E|=0$, there exists an outer function $f \in C^{1/2}(\T) \cap \mathcal{D}$ whose boundary zero set $\{ z \in \T : f(z) = 0 \}$ is exactly $E$. See \cite[Theorem 4.4.3]{dirichletspaceprimer}, or \cite{taylor1970ideals} for examples in $C^\infty(\T) \cap \mathcal{D}$. For the particular $E$ constructed in \thref{P:BScounterexample}, any such function is non-cyclic in $\mathcal{D}$. This should be contrasted with the construction of Brown and Cohn in \cite{browncohn1985} that exhibits a cyclic function in $C(\T) \cap \mathcal{D}$ which vanishes on any prescribed compact polar set $E \subset \T$. Thus the same set $E$ is the exact boundary zero set of two continuous outer Dirichlet functions, one cyclic and one non-cyclic. It follows that \thref{C:Noncyclicity} fails if $C^{\alpha}(\T)$ is replaced by the larger class $C(\T)$.

\subsection{The heart of the matter}

\subsubsection{Carleson integrals.}
For an arc $I\subseteq\T$ and a nonzero function $f\in\mathcal D$, consider the Carleson integral over $I$, namely
\begin{equation}\label{E:CarlDef}
 \Carl_I(f)=\int_I\int_\T
 \frac{(|f(t)|^2-|f(s)|^2)(\log |f(t)|^2-\log |f(s)|^2)}
 {|t-s|^2}\,d\m(s)\,d\m(t).
\end{equation}

According to \cite[Theorem 7.4.1]{dirichletspaceprimer}, for outer functions $f$ we have the identity $\Carl_\T(f) = 2\int_\D |f'|^2 \, dA$. We remark that we will use the Carleson integrals above even for functions $f \in \mathcal D$ which are not outer. Indeed, the estimate in \thref{P:MainTraceEstimate} below applies anyway. This will play a role, albeit minor.

\subsubsection{The main estimate.}

For a family $\mathcal{I}$ of intervals on $\T$ with disjoint interiors, we introduce a notation for the sum of averages of $|f|^2$ over the family:
\begin{equation} \label{E:TraceFuncDef}
\mathcal{T}_\mathcal{I}(f) = \sum_{I \in \mathcal{I}} \frac{1}{|I|}\int_{I} |f|^2 \, d\m.    
\end{equation}

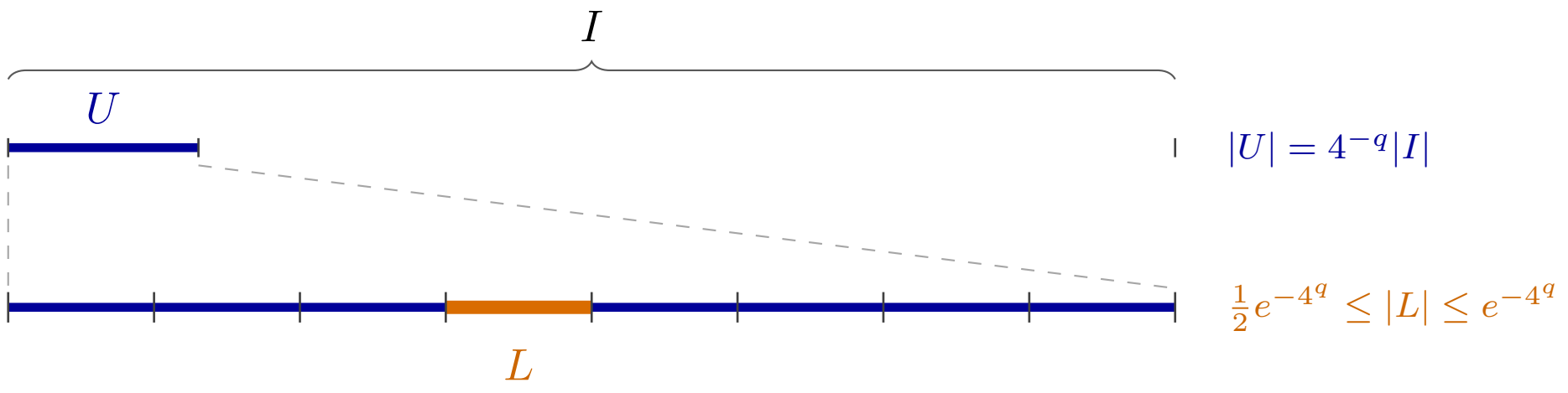
\begin{figure}[!ht]
\centering
\begin{tikzpicture}[
  x=1cm, y=1cm, line cap=butt,
  retained/.style={draw=blue!60!black, line width=2pt},
  gap/.style={draw=green!55!black, line width=.85pt, line cap=round,
              dash pattern=on 0pt off 2.3pt},
  guide/.style={draw=black!35, dashed, line width=.4pt},
  tick/.style={draw=black!75, line width=.5pt},
  length label/.style={anchor=west, font=\footnotesize, inner sep=0pt}
]
  \draw[decorate, decoration={brace,amplitude=4pt}, draw=black!65]
    (0,.54) -- (9.2,.54) node[midway,above=5pt,text=black] {$I$};
  \draw[retained] (0,0) -- (1.50,0);
  \foreach \x in {0,1.50,9.2}
    \draw[tick] (\x,-.075) -- (\x,.075);
  \node[above=2pt,text=blue!60!black] at (.75,0) {$U$};
  \node[length label,text=blue!60!black] at (9.60,0)
    {$|U|=4^{-q}|I|$};

  \draw[guide] (0,-.14) -- (0,-1.11);
  \draw[guide] (1.50,-.14) -- (9.2,-1.11);
  \draw[retained] (0,-1.26) -- (9.2,-1.26);
  \draw[orange!85!black,line width=3pt] (3.45,-1.26) -- (4.60,-1.26);
  \foreach \k in {0,...,8}
    \draw[tick] ({1.15*\k},-1.38) -- ({1.15*\k},-1.14);
  \node[below=2.5pt,text=orange!80!black] at (4.025,-1.38) {$L$};
  \node[length label,text=orange!80!black] at (9.60,-1.26)
    {$\tfrac12 e^{-4^q}\le |L|\le e^{-4^q}$};
\end{tikzpicture}
\caption{The lower row shows the subdivision of $U$.
In the construction of our Carleson set $E$, the same procedure will be recursively applied to each $L$.}
\label{fig:cantor-subdivision}
\end{figure}

The next proposition is a variant of \cite[Theorem 11.1]{chen2026brownshields}. It may be beneficial to keep referring to Figure~\ref{fig:cantor-subdivision} during the reading of the following parts in the red box. \footnote{The proofs of \thref{P:MainTraceEstimate} and \thref{C:MainTraceEstimateCor},
as well as the iteration in \eqref{E:IterationEq}, use only the boundary
modulus of the function and so do not require analyticity. Analyticity enters only when we bound
the Carleson integral by the Dirichlet integral. We mark this with a red exclamation mark below.}

\begin{tcolorbox}[
  colback=red!3,
  colframe=red!28,
  breakable
]
\begin{prop}\thlabel{P:MainTraceEstimate}
    Let $I \subseteq \T$ be an interval and let $q$ be a large positive integer. Let $U$ be the leftmost subinterval of $I$ of length $4^{-q}|I|$. Divide $U$ into equal-length intervals $L$ satisfying
    \begin{equation}
        \label{E:Llength}
        \frac{e^{-4^q}}{2} \leq |L| \leq e^{-4^q}.
        \end{equation} Let the finite family $\mathcal{I}'$ consist of all the intervals $L$ of this division. Then, for every $f \in C(\T) \cap \mathcal D$, we have
    \begin{equation} \label{E:traceEstimateLocal}
        \frac{1}{|I|}\int_I |f|^2\,d\m
        \leq 16 \pi^2 \cdot \Carl_I(f)
        + \frac{e^{-4^q/2}4^q}{|I|}\mathcal{T}_{\mathcal{I'}}(f).
    \end{equation}
\end{prop}

\begin{proof}
    Note first that for $a,b>0$ we have
    \begin{equation}\label{E:ScalarIneq}
       a\leq \frac{4}{4^q}(a-b)(\log a-\log b)+e^{4^q/2}b.
    \end{equation}
    Indeed, if $a\leq e^{4^q/2}b$, the inequality is clear, for the first term on the right-hand side is non-negative. Otherwise $a > e^{4^q/2}b$, and by taking logarithms and rearranging it follows that $\frac{2(\log a - \log b)}{4^q} > 1$. Since $e^{4^q/2} > 2$, the inequality $a > e^{4^q/2}b$ also implies $2(a-b) \geq a$. Then \eqref{E:ScalarIneq} follows by multiplying the two derived inequalities.

    In \eqref{E:ScalarIneq} take $a=|f(t)|^2$ and $b=|f(s)|^2$, and integrate over $t\in I$ and $s \in U$. We obtain

   \begin{equation}
        \label{E:ScalarIneqIntegrated}
       |U| \int_I |f|^2 d\m \leq \frac{4}{4^q} \int_I\int_U
 (|f(t)|^2-|f(s)|^2)(\log |f(t)|^2-\log |f(s)|^2)\,d\m(s)\,d\m(t) +  e^{4^q/2} |I|\int_U |f|^2 d\m.
  \end{equation}
    In the first integral term on the right-hand side, we multiply the integrand by $\frac{4\pi^2|I|^2}{|t-s|^2}$ which is greater than $1$ since $|t-s| < 2\pi|I|$ as $t,s \in I$, and then we extend the $s$-integral from $U$ to $\T$ by non-negativity of the integrand. For the second integral term, we use that $U$ is the union of the intervals $L \in \mathcal{I}'$ and \eqref{E:Llength} to obtain     
 
      \[  e^{4^q/2} |I|\int_U |f|^2 d\m = e^{4^q/2} |I| \sum_{L \in \mathcal{I}'} \int_{L} |f|^2 d\m \leq e^{-4^q/2} |I| \mathcal{T}_{\mathcal{I'}}(f).\]
      Putting these estimates into \eqref{E:ScalarIneqIntegrated} and dividing by $|I|\cdot|U|$ gives
      \begin{align*}\label{E:AverageEstimate}\frac{1}{|I|}\int_I |f|^2 \, d\m &\leq \frac{16 \pi^2 |I|}{4^q\cdot|U|}\Carl_I(f) + \frac{e^{-4^q/2}}{|U|}\mathcal{T}_{\mathcal{I'}}(f) \\
      &= 16 \pi^2 \Carl_I(f) + \frac{e^{-4^q/2}4^q}{|I|}\mathcal{T}_{\mathcal{I'}}(f).
      \end{align*}      
\end{proof}     
\end{tcolorbox}

\begin{cor} \thlabel{C:MainTraceEstimateCor}Let $\mathcal{I} = \{I_1, \ldots, I_s\}$ be a finite family of intervals in $\T$ of equal length with disjoint interiors. For each $I_j \in \mathcal{I}$, let $\mathcal{I}'_j(q)$ be the interval family constructed in \thref{P:MainTraceEstimate}, with common interval length across all those families, and let $\mathcal{I}'(q)$ be the union of those families. For all $q$ sufficiently large, it holds that for any $f\in C(\T)\cap\mathcal D$, we have
    \begin{equation}
        \label{E:traceEstimateGlobal}
        \mathcal{T}_{\mathcal{I}}(f) \leq 16 \pi^2 \Carl_\T(f) + \frac{1}{4}\mathcal{T}_{\mathcal{I}'(q)}(f).
    \end{equation} Given $\epsilon > 0$, for all sufficiently large $q$ we have that the capacity of $\bigcup_{L \in \mathcal{I'}(q)} L$ is less than $\epsilon$.
\end{cor}

\begin{proof}
     Choose $q$ at least so large that $\frac{e^{-4^q/2}4^q}{|I_j|} < 1/4$. The intervals $I_j$ have disjoint interiors, so $\sum_{j=1}^s\Carl_{I_j}(f)\leq\Carl_\T(f)$ by additivity of \eqref{E:CarlDef}. Summing \eqref{E:traceEstimateLocal} over $j$ gives \eqref{E:traceEstimateGlobal}.
    The family $\mathcal{I}_j'(q)$ consists of intervals contained within a common interval $U_j$ of length $4^{-q}|I_j|$. Thus the capacity of $\cup_{L \in \mathcal{I}'_j(q)} L$ tends to $0$ as $q \to \infty$. The capacity claim then follows, since there are only finitely many families $\mathcal{I}_j'(q)$.
\end{proof}

\subsection{The construction of the Carleson set $E$.}

\subsubsection{Carleson sets.}

For a finite or countable set of arcs $\mathcal{G}$, the Carleson entropy is
\begin{equation}
    \label{E:CarlesonEntropyDef}
    \kappa(\mathcal{G}) = \sum_{G \in \mathcal{G}} |G| \log(e/|G|).
\end{equation}

A compact set $E \subset \T$ with $|E|=0$ is a Carleson set if $\kappa(\mathcal{G}(E)) < \infty$, where $\mathcal{G}(E)$ is the family of maximal connected intervals constituting the complement $\T \setminus E$.

We have the \textbf{merging law}: if each interval in $\mathcal{G'}$ is a disjoint union of intervals in $\mathcal{G}$ (say, up to a set of measure zero, and each interval in $\mathcal{G}$ being used in merging only once), then $\kappa(\mathcal{G'}) \leq \kappa(\mathcal{G})$. This follows immediately from the inequality $\phi(a+b) \leq \phi(a) + \phi(b)$ for $\phi(t) = t\log(e/t)$. 

\subsubsection{Proof of \thref{P:BScounterexample}}

We construct $E$ as an intersection of unions of intervals in the families appearing in \thref{C:MainTraceEstimateCor}, and iterate the estimate in \eqref{E:traceEstimateGlobal} to prove \eqref{E:SmallAverageFromVanishing}. Finally, we verify that $E$ is Carleson. The construction is a variant of the one given in \cite[Section 11.3]{chen2026brownshields}.

Set $\mathcal I_1=\{\T\}$. At each step $n\geq2$, apply\footnote{At step $n=2$ we have $I = \T$ and so we have no endpoints. Take $U$ in \thref{P:MainTraceEstimate} equal to $\{ e^{i2\pi x} \in \T : x \in [0, 4^{-q_2}]\}$.} \thref{C:MainTraceEstimateCor} to $\mathcal I_{n-1} = \mathcal I$ to obtain $\mathcal I_n = \mathcal I'(q_n)$, choosing $q_n$ so that the capacity of $\bigcup_{L\in\mathcal I_n}L$ tends to $0$, and so that we have
\begin{equation}\label{E:IterationIneq}
        \mathcal{T}_{\mathcal{I}_{n-1}}(f) \leq 16 \pi^2 \Carl_\T(f) + \frac{1}{4}\mathcal{T}_{\mathcal{I}_{n}}(f)
\end{equation} 

If we take (as we may) the intervals $L$ to be closed, then
\[E = \bigcap_{n\geq 1} \Big( \bigcup_{L \in \mathcal{I}_n} L \Big)\] is clearly a closed polar set. The left endpoint of each $L \in \mathcal{I}_n$ survives to $E$. This follows immediately from the construction of the families $\mathcal{I}_n$.

Let us prove \eqref{E:SmallAverageFromVanishing} for $h$ as stated in the proposition. If a function $h \in C(\T) $ satisfies $|h(t)| \leq C_h \cdot \dist{t}{E}^{1/2}$ for some constant $C_h > 0$, then since each $L \in \mathcal{I}_n$ intersects $E$, we have $\dist{t}{E} \leq 2\pi |L|$ for $t \in L$, and consequently $|h(t)| \leq  C_h 2\pi \sqrt{|L|}$ for $t \in L$. According to \eqref{E:TraceFuncDef}, we obtain
    \begin{equation}
        \label{E:LocalHest}
        \mathcal{T}_{\mathcal{I}_n}(h) \leq \sum_{L \in \mathcal{I}_n} \frac{1}{|L|} \int_L C_h^2 4\pi^2 |L| \, d\m \leq C_h^2 4 \pi^2 \sum_{L \in \mathcal{I}_n} |L| \leq 4 \pi^2 C_h^2.
    \end{equation}
    The last inequality follows from disjointness of the interiors of the intervals $L$ and $|\T|=1$. Iterating \eqref{E:IterationIneq} and using \eqref{E:LocalHest} we obtain that for any $p \geq 2$ we have
\begin{equation}\label{E:IterationEq}
\begin{aligned}
    \int_\T |h|^2\,d\m
    &\leq 16\pi^2 \left(\sum_{j=0}^{p-2} 4^{-j}\right)
          \Carl_\T(h)
          + \frac{\mathcal{T}_{\mathcal{I}_p}(h)}{4^{p-1}} \\
    &\leq C\Carl_\T(h)
          + \frac{4\pi^2 C_h^2}{4^{p-1}}
\end{aligned} 
\end{equation}
with $C$ some (finite) numerical constant.
Sending $p\to\infty$ in \eqref{E:IterationEq} we obtain $\int_\T |h|^2\,dm  \leq C\Carl_\T(h)$. If $h = IQ$ is the inner-outer factorization of $h$ into inner factor $I$ and outer factor $Q$, then $|h| = |Q|$ almost everywhere on $\T$, and so $\int_\T |h|^2\,dm = \int_\T |Q|^2\,dm$ and $\Carl_\T(h) = \Carl_\T(Q)$. Hence 
\begin{align*}
    |h(0)|^2 \leq \int_\T |h|^2\,dm &\leq C\Carl_\T(Q)\\
    &\overset{\textcolor{red}{!}}{=} 2C\int_\D |Q'(z)|^2\,dA(z)\\
    &\leq 2C\int_\D |h'(z)|^2\,dA(z).
\end{align*}
The second line uses the remark that follows \eqref{E:CarlDef}, and the last line uses that division by an inner factor can only decrease the Dirichlet integral \cite[Corollary 7.6.3]{dirichletspaceprimer}. This proves \eqref{E:SmallAverageFromVanishing}.

It remains to verify that $E$ is indeed a Carleson set. For $I \in \mathcal{I}_{n}$ we have according to Figure~\ref{fig:cantor-subdivision} that $I \setminus \Big( \bigcup_{I' \in \mathcal{I}_{n+1}, I' \subset I} I'\Big) = I \setminus U$. Each interval in $\mathcal{G}(E)$ is, up to endpoints, a countable disjoint union of the intervals $I \setminus U$, $I \in \mathcal{I}_{n}$, $n\geq1$. By the merging law, it suffices to bound the sum of the Carleson entropies of the intervals $I \setminus U$. For every $I$, the Carleson entropy of $I \setminus U$ is
$|I \setminus U| \log(e/|I \setminus U|) \leq |I| \log(e/|I|)$.

Let $N_n$ be the number of intervals in $\mathcal I_n$, and let $\ell_n$ be their common length. By the above paragraph, each $I\in\mathcal I_n$ contributes not more than $\ell_n\log(e/\ell_n)$ to the entropy of $\mathcal G(E)$. In order to show that $E$ is a Carleson set, it thus remains to show that
\begin{equation}
    \label{E:EEntropyExpression}
    \sum_{n\geq1} N_n\ell_n\log(e/\ell_n)<\infty.
\end{equation}
For $n\geq2$, the ratio of $N_n\ell_n$ to $N_{n-1}\ell_{n-1}$ is the fraction of each parent interval $I\in\mathcal I_{n-1}$ covered by the blue interval $U$ in Figure~\ref{fig:cantor-subdivision}. Thus
\begin{equation}
    \label{E:GeometricDecrease}
    \frac{N_n\ell_n}{N_{n-1}\ell_{n-1}}= 4^{-q_n}.
\end{equation}
By \eqref{E:Llength} we have $\ell_n \geq e^{-4^{q_n}}/2$, so
\[
    \log(e/\ell_n) \leq 2 + 4^{q_n} \leq 2\cdot 4^{q_n}.
\]
Consequently,
\[
    N_n\ell_n\log(e/\ell_n)
    \leq 2 N_{n-1}\ell_{n-1}.
\]
But by \eqref{E:GeometricDecrease}, $N_{n}\ell_{n}$ is certainly at least geometrically decreasing in $n$. Thus the last displayed inequality shows that \eqref{E:EEntropyExpression} holds. The proof is complete.

\subsection*{Acknowledgements}
The author thanks Linus Bergqvist and Daniel Seco
for their helpful comments on the manuscript, and Adem Limani for discussions on the role of Carleson sets in Dirichlet cyclicity problems.


\begin{thebibliography}{EKMR14}
\bibitem[BC85]{browncohn1985}
L. Brown and W. Cohn, Some examples of cyclic vectors in the Dirichlet space, \emph{Proc. Amer. Math. Soc.} \textbf{95} (1985), no. 1, 42--46.
\bibitem[CFG26]{chen2026brownshields}
B. Chen, X. Fang and F. Guo, The Brown--Shields conjecture: Two approximation theorems and a counterexample, arXiv:2611.16997 (2026).
\bibitem[EKMR14]{dirichletspaceprimer}
O. El-Fallah, K. Kellay, J. Mashreghi and T. Ransford, \emph{A primer on the Dirichlet space}, Cambridge Tracts in Mathematics, vol. 203, Cambridge University Press, Cambridge, 2014.
\bibitem[TW70]{taylor1970ideals}
B. A. Taylor and D. L. Williams, Ideals in rings of analytic functions with smooth boundary values, \emph{Canad. J. Math.} \textbf{22} (1970), no. 6, 1266--1283.
\end{thebibliography}
\end{document}